\documentclass[a4paper,twoside,11pt]{amsart}
\newcommand{\Kurztitel}{Degree for Dimension}
\usepackage{fullpage}
\usepackage{amsmath} \usepackage{amssymb} \usepackage{amsopn}
\usepackage{amsthm} \usepackage{amsfonts} \usepackage{mathrsfs}

\usepackage[pdfpagelabels,pdftex]{hyperref}
\usepackage[dvips]{graphicx} \usepackage[usenames]{color}
\usepackage[all]{xy}

\hypersetup{
  pdftitle={},
  pdfauthor={},
  pdfsubject={},
  pdfkeywords={},
  colorlinks=false,
  breaklinks=true,
  bookmarksopen=true,
  bookmarksnumbered=true,
  pdfpagemode=UseOutlines,
  plainpages=false}

\DeclareMathOperator{\rH}{H}

\DeclareMathOperator{\rK}{K}

\DeclareMathOperator{\rR}{R}

\newcommand{\bQ}{{\mathbb Q}}
\newcommand{\bR}{{\mathbb R}}

\newcommand{\cZ}{{\mathscr Z}}

\newcommand{\surj}{\twoheadrightarrow}

\DeclareMathOperator{\id}{id}
\DeclareMathOperator{\pr}{pr}
\newcommand{\ev}{{\rm ev}}

\DeclareMathOperator{\Out}{Out}
\DeclareMathOperator{\Hom}{Hom}
\DeclareMathOperator{\Inn}{Inn}
\DeclareMathOperator{\Aut}{Aut}

\DeclareMathOperator{\Pic}{Pic}

\DeclareMathOperator{\Br}{Br}

\DeclareMathOperator{\Ext}{Ext}

\DeclareMathOperator{\res}{res}

\DeclareMathOperator{\ind}{ind}
\DeclareMathOperator{\Ind}{Ind}

\DeclareMathOperator{\Gal}{Gal}

\newcommand{\ph}{\varphi}
\newcommand{\teta}{\vartheta}

\newcommand{\sep}{{\rm sep}}
\newcommand{\alg}{{\rm alg}}
 
\newcommand{\sh}{{\rm sh}}

\newtheorem{thm}{Theorem}

\newtheorem{prop}[thm]{Proposition}
\newtheorem{lem}[thm]{Lemma}

\newtheorem{cor}[thm]{Corollary}
\newtheorem{conj}[thm]{Conjecture}

\theoremstyle{definition}

\theoremstyle{remark}
\newtheorem{rmk}[thm]{Remark}

\newenvironment{pro}[1][Proof]{{\it{#1:}} }{\hfill $\square$}
\newenvironment{pro*}[1][Proof]{{\it{#1:}} }{}
\newenvironment{pro**}[1][]{{\it{#1}} }{\hfill $\square$}

\newcounter{absatzcounter}[section]
\numberwithin{equation}{section}

\begin{document}
%




\hrule width\hsize
\hrule width\hsize
\hrule width\hsize

\vspace{1.5cm}

\title[\Kurztitel]{Trading degree for dimension \\[1ex] in the section conjecture: \\[1ex] The non-abelian Shapiro Lemma} 
\author{Jakob Stix}
\thanks{The author acknowledges support provided by DFG grant 
STI576/1-(1+2).}
\address{Mathematisches  Institut, Universit\"at Heidelberg, Im Neuenheimer Feld 288, 69120 Heidelberg}
\email{stix@mathi.uni-heidelberg.de}
\urladdr{http://www.mathi.uni-heidelberg.de/~stix/}
\date{February 10, 2009} 

\maketitle

\begin{quotation} 
  \noindent \small {\bf Abstract} --- This note aims at providing evidence for the section conjecture of anabelian geometry by establishing its behaviour under Weil restriction of scalars. In particular,  the \'etale fundamental group of the Weil restriction is determined by means of a Shapiro Lemma for non-abelian group cohomology.
\end{quotation}


\setcounter{tocdepth}{1} {\scriptsize \tableofcontents}


\section{Introduction and results}
Let $K$ be a field with fixed separable  closure $K^\sep$.
The \'etale fundamental group $\pi_1(X)$ of a geometrically connected variety $X/K$
forms naturally a non-abelian extension of  pro-finite groups
\[ 1 \to \pi_1(X \otimes K^\sep) \to \pi_1(X) \to \Gal(K^\sep/K) \to 1,\]
that  we abbreviate by $\pi_1(X/K)$.

\subsection{The section conjecture}
A $K$-rational point $x \in X(K)$ yields by functoriality a section $\Gal(K^\sep/K) \to \pi_1(X)$ of $\pi_1(X/K)$, with image the decomposition group of a point $\tilde{x}$ above $x$ in the universal pro-\'etale cover of $X$. Having neglected base points and due to the choice of $\tilde{x}$, only the class of a section up to conjugation by elements from $\pi_1(X \otimes K^\sep)$ is well defined. Let us denote by $S_{\pi_1{(X/K)}}$ the set of conjugacy classes of sections of $\pi_1(X/K)$. The section conjecture of Grothendieck's anabelian geometry \cite{letter}  speculates the following.
\begin{conj}[Grothendieck]
The map $X(K) \to S_{\pi_1(X/K)}$ which sends a rational point to the section given by its conjugacy class of decomposition groups is bijective if $K$ is a number field and $X/K$ is a geometrically connected, smooth, projective curve of genus at least $2$.
\end{conj} 
There is also a version of the section conjecture for affine curves. Here rational points at infinity will lead to an abundance of additional \textit{cuspidal} sections, see \cite{eh} and \cite{stix:cuspidalex}. But apart from the obvious modification forced upon us by cuspidal sections the conjecture remains the same. The condition on the genus gets replaced by asking the Euler-characteristic to be negative. 

A birational version of the section conjecture over $p$-adic local fields was successfully addressed by Koenigsmann in \cite{Koenigsmann}, and later by Pop \cite{pop:propsc} in a truncated version that exploits spectacularly modest pro-finite data to recover rational points.

Only recently evidence for the section conjecture could be found through the first examples of curves which satisfy the conjecture, see \cite{stix:periodindex} and \cite{hs}, though for the reason of having neither points nor sections. Another source of evidence has been the study of the cycle class of a section  as pioneered by \cite{ew}.

The goal of the present paper is to provide evidence for the section conjecture from a different direction. 

\subsection{Results --- trading degree for dimension} The evidence for the section conjecture presented in this note consists in its compliance with Weil restriction of scalars, see Section \ref{sec:Weil}. Of course, in order for this to make sense, we widen the applicability of the conjecture beyond the case of curves.

Let $L/K$ be a finite separable field extension within $K^\sep$, so that $\Gal_L=\Gal(K^\sep/L)$ is a subgroup of $\Gal_K=\Gal(K^\sep/K)$. Let $X/L$ be a quasi-projective, geometrically connected variety and $\rR_{L/K} X$ its Weil restriction of scalars as a geometrically connected variety over $K$. In Section \ref{sec:nabc} we will construct an induction functor for extensions which turns out to describe the fundamental group of the Weil restriction in characteristic $0$ as an extension as follows.

\begin{thm}[Theorem \ref{thm:pi1weiltext} in Section \ref{sec:Weil}]  \label{thm:pi1weil}
Let $K$ be a field of characterisitc $0$ or let $X/L$ be projective. Then the fundamental group $\pi_1(\rR_{L/K} X/K)$ of the Weil restriction $\rR_{L/K} X$ of scalars  is isomorphic to the non-abelian induction $\Ind_{\Gal_L}^{\Gal_K} \pi_1(X/L)$.
\end{thm}

Next, a non-abelian analogue of Shapiro's Lemma yields a description of the set of conjugacy classes of sections for an induction.

\begin{thm}[Corollary \ref{cor:sections} in Section \ref{sec:Weil}] \label{thm:corsections} 
Let $E = [1\to N \to E \to H \to 1]$ be an extension with  induction $\Ind_H^G(E)=[1\to M \to \Ind_H^G(E) \to G \to 1]$ with respect to a subgroup $H\subseteq G$. 

Then  $N$-conjugacy classes of sections of $E \surj H$ are naturally in bijection with $M$-conjugacy classes of sections of $\Ind_H^G(E) \surj G$.
\end{thm}

Combining Theorem \ref{thm:pi1weil} and Theorem \ref{thm:corsections} above, we obtain our piece of evidence for the section conjecture. We note in passing, that no special assumption on the geometry of the smooth, quasi-projective variety $X$ is used.
\begin{thm} \label{thm:scweil}
Let $L/K$ be a finite separable field extension and let $X/L$ be a quasi-projective, geometrically connected  variety. Let $K$ have characteristic $0$ or let $X/L$ be projective.  Then applying the functor $\pi_1$ yields a bijective map $X(L) \to S_{\pi_1{(X/L)}}$ if and only if it yields a bijective map $\rR_{L/K} X(K) \to S_{\pi_1{(\rR_{L/K} X/K)}}$. 
\end{thm}

\begin{cor}
The section conjecture holds for smooth projective curves over number fields if it holds for smooth, projective algebraic $\rK(\pi,1)$ spaces over $\bQ$ (see \cite{stix:diss} Appendix A), which embed into their Albanese variety and have  non-vanishing Euler-Poincar\'e characteristic. 
\end{cor}

This corollary explains the title of the article. We have lowered the degree of the number field in the section conjecture to  $1$ at the expense of working with varieties of dimension exceeding~$1$.  Of course, these higher dimensional varieties of interest are simply $\bQ$-forms of products of smooth, projective hyperbolic curves and so the trade might be marginal. But on the one hand, we are not required to limit the section conjecture to curves or products of curves, and secondly, the arithmetic of $\bQ$ and so presumably also sections over $\Gal_\bQ$ are arithmetically much simpler than for more general algebraic number fields. For exampe, the recent modularity results were first proven for representations of the full $\Gal_\bQ$. So the base has become simpler. 

Another argument in favour of the improvement of our situation is the following weak analogue of 
\cite{stix:periodindex} Theorem 17.
\begin{prop}
Let $X/\bQ$ be a smooth, projective and geometrically connected variety which is an algebraic $\rK(\pi,1)$ space such that $\pi_1(X/\bQ)$ admits a section. Then any Galois invariant line bundle $L \in \Pic_X^0(\bQ)$ has vanishing Brauer obstruction $b(L)=0 \in \Br(\bQ)$, hence belongs to a genuine line bundle on $X$.
\end{prop}
\begin{pro}
The local components of $b(L)$ are obtained by base change to $\bQ_p$ or $\bR$ which preserves the assumptions. The real section conjecture as in \cite{stix:periodindex} Theorem 24 applies also for higher dimensional $\rK(\pi,1)$ spaces and yields the existence of a real point, and so $b(L)_\bR=0$. 

Over a $p$-adic field the argument of \cite{stix:periodindex} Proposition 12 and Corollary 14 still show  that the order of $b(L)_{\bQ_p}$ is a power of $p$. 

The result now follows from the global reciprocity which states that the local invariants of $b(L)$ sum up to $0$. A sum of summands of prime power order where every prime occurs at most once can only vanish if all the summands vanish. Hence $b(L)$ vanishes by the local global principle for Brauer groups of number fields.
\end{pro}


\section{Non-abelian cohomology} \label{sec:nabc}

\subsection{Twisted generalized wreath products} \label{subsec:wreath}
We recall some group theory in order to fix notations and to put it into a form useful for the sequel. 
Although we will ultimately apply the results in the context of pro-finite groups, we choose to neglect the topology in the presentation, because this frees us from adding an abundance of "continuous" everywhere. However we note, that in the pro-finite case the subgroups in question should be closed subgroups.

\subsubsection{Wreath products} The \textbf{wreath product} of two groups $G$ and $N$ along a right $G$-set $A$ is the semi\-direct product 
\[N \wr G := (\prod_{\alpha \in A} N) \rtimes G\]
 with respect to the action of $G$ on $\prod_{\alpha \in A} N$ given by $g.(n_\alpha)_{\alpha\in A} = (n_{\alpha g})_{\alpha \in A}$, see \cite{huppert1} \S15.6.

\subsubsection{Induction for groups acting on groups} Let $G$ be a group. A \textbf{$G$-group} is a group $N$ together with an action $\teta: G \to \Aut(N)$. For a subgroup $H$ in $G$, induction is a functor $\ind_H^G$ from $H$-groups to $G$-groups defined as follows. As a group 
\[\ind_H^G(N) =  \{f:G \to N; \ f(hg) = \teta(h)(f(g)) \text{ all $h \in H$} \}\]
 with pointwise multiplication. The $G$-action on $f\in \ind_H^G(N)$ comes from right translation of the argument, so $(g.f)(\alpha) = f(\alpha g)$. 

\subsubsection{Twisted generalized wreath products} The \textbf{twisted generalized wreath product} of the group $G$ with subgroup $H$ relative to the $H$-group $N$ is the semidirect product 
\[N \wr_H G := (N,\teta) \wr_H G := \ind_H^G(N) \rtimes G\]
with respect to the natural $G$-action of the induction as a $G$-group, see \cite{neumann} \S2, \cite{meldrum} I.9.3, \cite{huppert1} \S15.10 and \cite{haran} \S1.

\subsubsection{Sections and \texorpdfstring{$\rH^1$}{first cohomology}} Let $N$ be a $G$-group. A $1$-cocycle of $G$ with values in $N$ is a map $a:G \to N$ such that for all $s,t \in G$ we have 
\[a_{st} = a_s (s.a_t).\]
 The first non-abelian cohomology $\rH^1(G,N)$ is the set of equivalence classes of $1$-cocycles, where cocycles $a$ and $b$ are equivalent if there is $c \in N$ with $a_s = c b_s (s.c)^{-1}$ for all $s \in G$. 
 
 Two sections of $N \rtimes G \surj G$ are equivalent if they differ by conjugation with an element of $N$.  The following lemma is well known and straight forward.
 
 \begin{lem} \label{lem:h1sect}
 The map which send a $1$-cocycle $s \mapsto a_s$ to the section $s \mapsto a_s \cdot s$, where we have identified $G$ with the second factor in $N \rtimes G$ establishes a natural bijection of 
 $\rH^1(G,N)$ with the set of equivalence classes of sections of $N \rtimes G \surj G$.
  \end{lem}

\subsubsection{The non-abelian Shapiro Lemma in degree 1} \label{subsec:sh1}
Let $N$ be an $H$-group for a subgroup $H \subset G$. The restriction of the $G$-group $\ind_H^G(N)$ to an $H$-group admits an $H$-equivariant map 
\[\ev_1:\ind_H^G(N)|_H \to N\]
 by evaluating at $1$.  The composition of restriction and evaluation at $1$ defines the Shapiro map 
 \[ \sh^1: \rH^1(G,\ind_H^G(N)) \to \rH^1(H,N).\]
 
 \begin{prop} \label{prop:shapiro1}
The Shapiro map $\sh^1: \rH^1(G,\ind_H^G(N)) \to \rH^1(H,N)$ is bijective. 
 \end{prop}
\begin{pro} A $1$-cocycle $s \mapsto b_s$ for $G$ with values in $\ind_H^G(N)$ is given by $b_{s,t}=b_s(t) \in N$ for all $s,t \in G$ such that (i) $b_{s,ht} = \teta(h)(b_{s,t})$ for all $s,t \in G$ and $h \in H$ and (ii) $b_{st,g} = b_{s,g}b_{t,gs}$ for all $s,t,g \in G$.
The map $\sh^1$ on the level of cocycles maps $b$ to $h \mapsto b_{h,1}$.

\textit{Surjectivity.} 
We choose a set of representatives $Y \subset G$ for $H\backslash G$ and obtain maps $\gamma : G \to H$ and $y:G \to Y$ such that $g=\gamma_gy_g$ for all $g\in G$. Let $a:H \to N$ be a $1$-cocycle, in particular $a_1=1$. We set
\[
b_{s,t} := \big(a_{\gamma_t} \teta(\gamma_t)(a_{y_t}) \big)^{-1}\big(a_{\gamma_{ts}} \teta(\gamma_{ts})(a_{y_{ts}}) \big)
\]
and a routine calculation shows that $b$ is a $1$-cocycle maping to $a$. The cocycle condition (ii) is best checked by noting that our definition of $b$ implies $b_{s,t} = (b_{t,1})^{-1}b_{ts,1}$.

\textit{Injectivity.}
Let $b,b'$ be cocycles with $\sh^1(b) \sim \sh^1(b')$. We have $f\in \ind_H^G(N)$ defined by $f(s) := (b'_{s,1})^{-1}b_{s,1}$ for all $s \in G$. It follows that $b'_{s,t} = f(t)b_{s,t}f(ts)^{-1}$ for all $s,t \in G$ which translates into $b\sim b'$.
\end{pro}

\begin{rmk} (1) An alternative proof is given in \cite{holt} Theorem 4 using the interpretation as conjugacy classes of  complements in semidirect products.

(2) Proposition \ref{prop:shapiro1} speaks about conjugacy classes of sections of a twisted wreath product $N \wr_H G = \ind_H^G(N) \rtimes G $. The assertion appears in \cite{pq} Thm 2.6 in the case of wreath products, i.e., trivial action of $H$ on $N$. The introuction of \cite{pq} contains the observation that this is a non-abelian version of Shapiro's Lemma but does not elaborate on this idea further.
\end{rmk}


\subsection{Extensions after Eilenberg and MacLane} We recall the theory of non-abelian extensions of Eilenberg and MacLane from \cite{emcl}.

\subsubsection{Kernel} A \textbf{kernel} or more precisely a \textbf{$G$-kernel} is a group $N$ together with an exterior action $\rho:G \to \Out(N)$ by a group $G$, where $\Out(N)=\Aut(N)/\Inn(N)$ is the group of exterior automorphisms of $N$. We denote the set of all $G$-kernels on $N$ by $\rK(G,N) = \Hom(G,\Out(N))$.
 
 \subsubsection{Center} Let $\rho:G \to \Out(N)$ be a $G$-kernel. The \textbf{center} $Z$ of the kernel  is the center of $N$ together with its inherited $G$-action $\chi=\chi_\rho:G \to \Aut(Z)$. For distinction purposes we may denote by $Z(\chi)$ the $G$-module $Z$ with module structure given by $\chi:G \to \Aut(Z)$. The set of kernels $\rho$ with center equal to $\chi$ is denoted by $\rK(G,N)_\chi$.
 
\subsubsection{Extensions} An \textbf{extension} $E$ of a group $G$ by a group $N$ is a short exact sequence \[1 \to N \to E \to G \to 1.\]
Isomorphisms of extensions of $G$ by $N$ respect both $G$ and $N$ identically. We denote the set of isomorphism classes of extensions by $\Ext(G,N)$. An extension $E$ leads to a kernel $\rho$ via the restriction to $N$ of the conjugation by preimages:
\[
\xymatrix@C+3ex@R-3ex{E \ar@{->>}[d] \ar[r]^(0.45){e(\ )e^{-1}|_N} & \Aut(N) \ar[d]\\
G \ar[r]^(0.4)\rho & \Out(N)}
\]
The set of isomorphy classes of extensions, whose kernel has center $\chi$ is denoted by $\Ext(G,N)_\chi$. The map $\Ext(G,N)_\chi \to \rK(G,N)_\chi$ that assigns to each extension its kernel is well defined.
 
 \subsubsection{Obstruction theory} A kernel which is the kernel of an extension is called \textbf{extendible}. By pullback of 
 \[ 1 \to N/Z \to \Aut(N) \to \Out(N) \to 1\]
 a kernel $\rho : G \to \Out(N)$ determines an extension 
 \[ 1 \to N/Z \to E_\rho \to G \to 1\]
 such that a presumptive extension $E$ with kernel $\rho$ by conjugation canonically sits in a diagram
 \[
 \xymatrix@R-2ex{ 1 \ar[r] & N \ar[d] \ar[r] & E \ar[r] \ar[d] & G \ar[r] \ar@{=}[d] & 1 \\
  1 \ar[r] & N/Z \ar[r] & E_\rho \ar[r] & G \ar[r] & 1}
 \]
 Let $\chi$ be the center of $\rho$. If $N$ were an abelian group with a $G$- submodule $Z$, then the existence of $E$ lifting $E_\rho$ were controlled by the coboundary $\delta(E_\rho)$ under 
 $\delta: \rH^2(G,N/Z) \to \rH^3(G,Z(\chi))$. The set of such lifts would receive a transitive action of the group $\rH^2(G,Z(\chi))$ via the homomorphism $\rH^2(G,Z(\chi)) \to \rH^2(G,N)$. In the non-abelian case discussed here (see \cite{emcl} \S7+8), as $Z$ is central in $N$, the same formulas with inhomogeneous cocycles which prove the assertions in the abelian case succeed to give the following result, except for the $0$ on both sides which follows from \cite{emcl} \S9+11.
 
\begin{prop}[essentially \cite{emcl}] \label{prop:ext}
The following is exact
\[
0 \to \rH^2(G,Z(\chi)) \to \Ext(G,N)_\chi \to \rK(G,N)_\chi \xrightarrow{\delta} \rH^3(G,Z(\chi)) \to 0
\]
in the sense that $\rH^2(G,Z(\chi))$ acts freely  on $\Ext(G,N)_\chi$ with quotient set equal to the set of extendible kernels $\delta^{-1}(0)$ and $\delta$ is surjective.
\end{prop}

The action of $a \in \rH^2(G,Z(\chi))$ on $\Ext(G,N)_\chi$ can be constructed on extensions as follows. Let $1 \to Z \to \cZ_a \to G \to 1$ be an extension realizing the cohomology class $a \in \rH^2(G,Z(\chi)),$ and let $1\to N \to E \to G \to 1$ be an extension with center of its kernel equal to $\chi$. Then $a.E$ equals the isomorphism class of the extension 
\[
1 \to N \xrightarrow{{\rm inclusion},0} (E \times_G \cZ_a)/\Delta(Z) \to G \to 1,
\]
where $\Delta: Z \to N \times Z$ is the antidiagonal $z \mapsto (z,-z)$.

\subsubsection{Categories of extensions}
Let $G$ be a group. The category $\Ext[G]$ has as objects extensions of $G$ with arbitrary kernel  and morphisms are maps of extensions up to composition by inner automorphisms from elements of the kernel. 

Pushing an extension $1 \to N \to E \to G \to 1$ by an automorphism of $N$ determines an action of $\Out(N)$ on $\Ext(G,N)$ such that the set of orbits equals the set $\Ext[G,N]$ of isomorphism classes in the category $\Ext[G]$ of extensions of $G$ by $N$. The map $\Ext(G,N) \to \rK(G,N)$ becomes $\Out(N)$-equivariant when $\Out(N)$ acts on $\rK(G,N)=\Hom(G,\Out(N))$ by composition with inner automorphisms of $\Out(N)$.


\subsection{Wreath product type extensions} \label{subsec:sh2}
 In this section we built on the work of Holt \cite{holt}.

\subsubsection{Wreath kernels}  
Let $N$ be a group and $H \subseteq G$ a subgroup. On 
\[M=\ind_H^G(N,1) = \prod_{\alpha \in H\backslash G} N\] 
we have an action of $G$ by 
\[
g.((n_\alpha)_{\alpha \in H\backslash G}) = (n_{\alpha g})_{\alpha \in H \backslash G}
\]
and an outer action $\ind_H^G(\Out(N),1) = \prod_{\alpha \in H\backslash G} \Out(N) \to \Out(M)$ given by 
\[
(f_\alpha)_{\alpha \in H \backslash G} (n_\alpha)_{\alpha \in H \backslash G} = (f_\alpha(n_\alpha))_{\alpha \in H \backslash G}.
\]
The two actions are compatible as follows. For $g \in G$ and $(f_\alpha)_{\alpha \in H \backslash G}  \in 
\prod_{\alpha \in H\backslash G} \Out(N) $ we have
\[
\big(g \cdot (f_\alpha)_{\alpha \in H \backslash G} \cdot g^{-1} \big) \big((n_\alpha)_{\alpha \in H \backslash G}\big) = \big(f_{\alpha g}(n_\alpha)\big)_{\alpha \in H \backslash G} = \big(g.(f_\alpha)_{\alpha \in H \backslash G}\big)\big((n_\alpha)_{\alpha \in H \backslash G}\big)
\]
resulting in a homomorphism
\[
R: (\Out(N),1) \wr_H G \to \Out(M).
\]

A \textbf{wreath kernel} of $G$ on $N$ is a kernel $\rho:G \to \Out(M)$ together with a lift along $R$ to a homomorphism $\tilde{\rho} : G \to (\Out(N),1) \wr_H G$ which is a section of the projection to $G$. Such a lift is unique if it exists, as two lifts differ at most by elements in $\prod_{\alpha \in H\backslash G} \Out(N)$ which injects into $\Out(M)$. The set of wreath kernels is thus a subset $\rK_{\rm wreath}(G,H; N) \subseteq \rK(G,M)$ of the set of all $G$-kernels on $M$.

\subsubsection{The Center of a wreath kernel}
 The center of $M$ equals $\ind_H^G(Z,1) = \prod_{\alpha \in H \backslash G} Z$ where $Z$ is the center of $N$. Thus the center of a wreath kernel lifts to a homomorphism 
 \[G \to (\Aut(Z),1) \wr_H G,\]
 which is a section of the projection to $G$.
The construction in Section \ref{subsec:sh1} performed on the level of cocycles, namely restriction to $H$ and then evaluation at $1$,  yields the two surjective maps $\sh$ as in the following diagram. 
\begin{equation} \label{eq:sh}
\xymatrix@R-2ex{\rK_{\rm wreath}(G,H; N) \ar[d]_{\text{center}} \ar@{->>}[r]^\sh & \rK(H,N) \ar[d]^{\text{center}} \\
\{\text{Sections of }  (\Aut(Z),1) \wr_H G \surj G\}\ar@{->>}[r]^(0.63)\sh & \Hom(H,\Aut(Z)) }
\end{equation}
The vertical maps associate to a kernel its center. In particular, the center of $M$ as a $G$-module under the center of the wreath kernel $\rho$ is nothing but $\ind_H^G(\chi)=\ind_H^G(Z(\chi))$ where $\chi$ is the center of the $H$-kernel $\sh(\rho)$.

\subsubsection{Wreath product type extensions}
A \textbf{wreath product type extension}, see \cite{holt} p.464, is an extension  of $G$ by $M=  \prod_{\alpha \in H\backslash G} N$, the kernel of which is a wreath kernel as above. The set of isomorphism classes of wreath product type extensions
 is denoted by $\Ext_{\rm wreath}(G,H;N)$.  We denote by $\Ext_{\rm wreath}(G,H;N)_\chi$ (resp.\ by $\rK_{\rm wreath}(G,H;N)_\chi$) the set of those  wreath type extensions whose kernel maps (resp.\ those wreath kernels which map) under $\sh$ to a kernel with center $\chi$.

\subsubsection{The Shapiro map for extensions} 
 Let $E = [1 \to M \to E \xrightarrow{\pr} G \to 1]$ be a wreath product type extension. The kernel of the map $\ev_1 : M \surj N$, which evaluates at $1$ is a normal subgroup of $E|_H = \pr^{-1}(H)$. We may therefore push the restriction of $E$ to $H$ by the map $\ev_1$ to obtain an extension of $H$ by $N$, that will be denoted $\sh^2(E)$.

\subsubsection{Non-abelian Shapiro Lemma in degree 2}

\begin{prop}
Let $H$ be a subgroup of $G$. Let $N$ be a group with center $Z$ and $H$-action $\chi$. 
We have a commutative ladder with exact rows in the sense as in Proposition \ref{prop:ext}
\[
\xymatrix@R-2ex@C-2.5ex@M+0.5ex{
0 \ar[r] & \rH^2(G,\ind_H^G(\chi)) \ar[r] \ar[d]^{\sh^2}_\cong &  \Ext(G,H;N)_\chi \ar@{->>}[d]^{\sh^2} \ar[r] & \rK_{\rm wreath}(G,H;N)_\chi \ar[r]^\delta \ar@{->>}[d]^{\sh} & \rH^3(G,\ind_H^G(\chi)) \ar[d]^{\sh^3}_\cong \ar[r] & 0 \\
0 \ar[r] & \rH^2(H,Z(\chi)) \ar[r] &  \Ext(H,N)_\chi \ar[r] & \rK(H,N)_\chi \ar[r]^\delta & \rH^3(H,Z(\chi)) \ar[r] & 0}
\]
where the vertical maps are induced by the respective Shapiro map, are all surjective and the two extremal ones are isomorphisms. 
\end{prop}
\begin{pro}
Exactness of the bottom row is Proposition \ref{prop:ext}. The commutativity of the diagram follows by tedious but elementary calculations on cochains. The abelian Shapiro Lemma shows that the two extremal vertical maps are isomorphisms. The surjectivity of $\sh$ was discussed in (\ref{eq:sh}).
Exactness of the top row follows again from Proposition \ref{prop:ext} besides  the surjectivity of $\delta$ which follows from a diagram chase. Now the surjectivity of the remaining vertical map follows again by diagram chase.
\end{pro}

\subsubsection{Split extensions}
The group $\prod_{\alpha \in H \backslash G} \Out(N)$ acts naturally by pushing on isomorphism classes $\Ext_{\rm wreath}(G,H;N)$ of wreath product type extensions. Let $\Ext_{\rm wreath}[G,H;N]$ denote the set of orbits. 
\begin{thm}[\cite{holt} Thm 3] \label{thm:holt}
The map $\sh^2$ yields a bijection $\Ext_{\rm wreath}[G,H;N] \xrightarrow{\sim} \Ext[H,N]$.
\end{thm}
This theorem by Holt has the following immediate corollary.
\begin{cor} \label{cor:existsplitting}
Let $E$ be a wreath product type extension in $\Ext_{\rm wreath}(G,H;N)_\chi$. Then $E$ splits as an extension of $E$ by $\prod_{\alpha \in H \backslash G} N$ if and only if $\sh^2(E)$ splits as an extension of $H$ by $N$.
\end{cor}


\subsection{Non-abelian induction of extensions} \label{subsec:nonabind}

\subsubsection{Adjoint to restriction}
Let $H$ be a subgroup of $G$. We are going to construct a right adjoint to the functor restriction from $G$ to $H$
\[
\res_H^G : \Ext[G] \to \Ext[H], \quad E \mapsto \res_H^G = E|_H
\]
which by analogy is called induction from $H$ to $G$ and denoted by $\Ind_H^G$. By Yoneda, all we have to achieve is a proof of  the following theorem.

\begin{thm} \label{thm:ind}
Let $F \in \Ext[H]$ be an extension of $H$.
The functor which sends $T \in \Ext[G]$ to $\Hom_{\Ext[H]}\big(\res_H^G(T), F\big)$ is representable. 
\end{thm}
A representing object as in the theorem is denoted by $\Ind_H^G(F)$. By definition,  adjointness 
\begin{equation} \label{eq:adj}
\Hom_{\Ext[H]}\big(\res_H^G(T),F\big) = \Hom_{\Ext[G]}\big(T,\Ind_H^G(F)\big)
\end{equation}
holds naturally, and $F \mapsto \Ind_H^G(F)$ is the sought for non-abelian induction functor.

If we apply adjointness to the extension $1 \to 1 \to G \xrightarrow{\id} G \to 1$ we get  the following immediate corollary.
\begin{cor} \label{cor:sections}
Let $H$ be a subgroup of $G$ and let $E = [1\to N \to E \to H \to 1]$ be an extension in $\Ext[H]$ with induction $\Ind_H^G(E)=[1\to M \to \Ind_H^G(E) \to G \to 1]$ in $\Ext[G]$. 
\begin{itemize}
\item[(1)] The adjointness map 
(\ref{eq:adj}) describes a bijection between $N$-conjugacy classes of sections of $E \surj H$ and $M$-conjugacy classes of sections of $\Ind_H^G(E) \surj G$.
\item[(2)] In particular, $E$ splits if and only if $\Ind_H^G(E)$ splits.
\end{itemize}
\end{cor}

\subsubsection{Proof of adjointness}  
\begin{pro**}
Let $F \in \Ext[H]$ be an extension of $H$. By Theorem \ref{thm:holt} there exist an extension $E \in \Ext[G]$ which is a wreath product type extension such that $\sh^2(E) \cong F$. We are done if we can find a bijection 
\begin{equation} \label{eq:adj2}
\sh:  \Hom_{\Ext[G]}\big(T,E\big) \xrightarrow{\sim} \Hom_{\Ext[H]}\big(\res_H^G(T),\sh^2(E)\big) 
\end{equation}
which is natural in $T = [1 \to S \to T \to G \to 1]$. 

Let $\ph:T \to E$ be a morphism of extensions. We define the map $\sh(\ph)$ as the composition $\ev_1 \circ \res_H^G(\ph)$ which is the vertical map in the following diagram.
\[
\xymatrix@R-2ex{
\res_H^G(T) \ar[d]^{\res_H^G(\ph)} & = & [1 \ar[r] & S \ar[r]  \ar[d] & T|_H \ar[d] \ar[r] & H \ar@{=}[d] \ar[r] & 1]\\ 
\res_H^G(E) \ar[d]^{\ev_1} & = & [1 \ar[r] &  M \ar[r]  \ar[d] & E|_H \ar[d] \ar[r] & H \ar@{=}[d] \ar[r] & 1]\\ 
\sh^2(E) & = & [1 \ar[r] &  N \ar[r]  & \sh^2(E) \ar[r] & H \ar[r] & 1]}
\]
We denote the projection $T \to G$ by $\pr$ and its restriction to $H$ by $\pr|_H : T|_H \to H$. A map $\ph: T \to E$ is nothing but an $M$-conjugacy class of sections of 
\[
\pr_1^\ast E = [ 1 \to M \to E \times_G T \to T \to 1],
\]
which is a wreath product type extension for the subgroup $T|_H \subseteq T$. On the other hand, a map $T|_H \to \sh^2(E)$ is an $N$-conjugacy class of sections of 
\[
(\pr|_H)^\ast(\sh^2(E))= [ 1 \to N \to \sh^2(E) \times_H T|_H \to H \to 1].
\]
Because of the isomorphism $\sh^2(\pr_1^\ast E) \cong (\pr|_H)^\ast(\sh^2(E))$, we get that in fact Theorem \ref{thm:ind} is equivalent to  Corollary \ref{cor:sections} with $F=\Ind_H^G(E)$ and $E=\sh^2(F)$, which is what we are going to prove in the sequel. Part (2) has already been shown in Corollary \ref{cor:existsplitting}. We may therefore assume that the wreath product type extension $F$ splits. 

We fix a splitting $\sigma$, which allows to lift the kernel $\rho_F : G \to \Out(M)$ to a homomorphism $\teta_F: G \to \Aut(M)$. By means of $\teta_F$ the extension $F$ is isomorphic to the semidirect product $M \rtimes_{\teta_F} G$. Because the following diagram is a fibre product diagram
\[
\xymatrix{ \big(\Aut(N),1\big) \wr_H G \ar[r] \ar[d]^{\tilde{R}} &  \big(\Out(N),1\big) \wr_H G  \ar[d]^{R} \\
\Aut(M) \ar[r] & \Out(M)}
\]
where the map $\tilde{R}$ is constructed analogously to the map $R$, the wreath kernel $\tilde{\rho}$, which lifts $\rho_F$, also lifts to a homomorphism $\tilde{\teta} : G \to \big(\Aut(N),1\big) \wr_H G$. 

Restriction of the section $\sigma$ to $H$ and evaluating at $1$ induce the corresponding section of $\sh^2(F)$ and a true action $\teta: H \to \Aut(N)$ by means of which $\sh^2(F) \cong N \rtimes_\teta H$. But $\teta$ also equals the evaluation at $1$ of the restriction to $H$ of $\tilde{\teta}$. Hence $M$ as a $G$-group via $\tilde{\teta}$ is nothing but $\ind_H^G(N)$ for the $H$-group $N$ via $\teta$. The result now follows from Proposition \ref{prop:shapiro1}.
\end{pro**}


\subsection{The anabelian case}
Although still a mystery in general, it is widely believed that a group only qualifies to be anabelian when  its center is trivial. In this section we will work out under the assumption, that the center $Z$ of $N$ is trivial, how the content of the Sections \ref{subsec:wreath} - \ref{subsec:nonabind} specialises.
Proposition \ref{prop:ext} yields a bijection
\[ \Ext(G,N) \to \rK(G,N)
\]
with inverse assigning to a kernel $\rho:G \to \Out(N)$ the pullback under $\rho$ of the extension $1 \to N \to \Aut(N) \to \Out(N) \to 1$. Hence all kernels are extendible and each extension is determined by its kernel up to isomorphism. 

It follows that an $N$-conjugacy class of sections of an extension with kernel $\rho:G \to \Out(N)$ is canonically the same as a lift of $\rho$ to an actual action $\tilde{\rho}:G \to \Aut(N)$ up to uniform conjugation by elements of $N$.

With $N$ also $M = \prod_{\alpha \in H \backslash G} N$ has trivial center. The fact that sections of $\big(\Out(N),1\big) \wr_H G \surj G$ and  sections of $\Out(N) \times H \surj H$ correspond to each other under $\sh^1$ up to conjugation by elements from $M$ (resp. $N$) thus explains Theorem \ref{thm:holt} in this case.
The content of Corollary \ref{cor:sections} follows by applying the same argument to the semi-direct products $\big(\Aut(N),1\big) \wr_H G \surj G$ and $\Aut(N) \times H \surj H$.


\section{Weil restriction of scalars} \label{sec:Weil}
We content ourselves with a discussion of Weil restriction of scalars relative a finite separable field extension for quasi-projective varieties.

\subsection{Properties of the Weil restriction}
Let $L/K$ be a finite separable field extension. Weil restriction of scalars relative $L/K$ is a functor 
$\rR_{L/K}$ from quasi-projective varieties over $L$ to quasi-projective varieties over $K$ that is right adjoint to scalar extension $-\otimes_K L$ and thus defined by an identification
\[
\Hom_L\big(Y \otimes_K L, X\big) = \Hom_K\big(Y,\rR_{L/K}(X)\big)
\] 
which is natural for schemes $X$ (resp. $Y$) quasi-projective over $L$ (resp. $K$), see \cite{BLR} VII.6. If $X$ has dimension $d$ then $\rR_{L/K}(X)$ has dimension $d[L:K]$.

Let $K'/K$ another field extension and $L'=L\otimes_K K' = \prod_\alpha L'_\alpha$ the decomposition of the tensor product in separable $K'$ extensions.
Then the following Mackey-formula holds
\[
\big(\rR_{L/K}(X)\big) \otimes_K K' = \prod_\alpha \rR_{L'_\alpha/K'}\big(X \otimes_L L'_\alpha\big).
\]
In particular, with an algebraic closure $K^\alg$ of $K$ we have 
\[
\big(\rR_{L/K}(X)\big) \otimes_K K^\alg = \prod_\sigma X \otimes_{L,\sigma} K^\alg,
\]
where the product is over all $K$-embeddings $\sigma: L \to K^\alg$. 
If $K^\alg$ contains $L$ a priori, then we have a prefered $\sigma_1 = \id$ which gives a projection 
\[
\pr_1 : \big(\rR_{L/K}(X)\big) \otimes_K K^\alg \to X \otimes_L K^\alg
\]
onto the corresponding  factor.

\subsection{The fundamental group of the Weil restriction} From now on we work with fields of characteristic $0$ or ask the variety $X/L$ to be projective.
In the extension $\pi_1\big(\rR_{L/K}(X)/K\big)$ the outer action by conjugation with lifts of $g \in \Gal_K$ on 
\[ \pi_1\big(\rR_{L/K}(X) \otimes_K K^\alg\big) = \prod_\sigma \pi_1\big(X \otimes_{L,\sigma} K^\alg\big)
\]
acts as $\pi_1(1 \otimes g^{-1})$ and therefore permutes the factors by mapping $\sigma \in \Hom_K(L,K^\alg) = \Gal_K/\Gal_L$ to $g^{-1}\sigma$. Reindexing the product with $\alpha = \sigma^{-1}$ transforms this to the right translation action on $\Gal_L \backslash \Gal_K$.  Thus we have established the following proposition.

\begin{prop}
The extension $\pi_1\big(\rR_{L/K}(X)/K\big)$ is a wreath product type extension with respect to the subgroup $\Gal_L \subseteq \Gal_K$.
\end{prop}

In the identification of the index set of the product with $\Gal_L \backslash \Gal_K$ we have chosen a distinguished embedding of $L$ in $K^\alg$. Evaluation at $1$ then is nothing but the $\pi_1(\pr_1)$ for 
\[
\pr_1 : \big(\rR_{L/K}(X)\big) \otimes_K L \to X 
\]
which is the adjoint map for the pair of adjoint functors $\rR_{L/K}$ and $-\otimes_K L$ and an $L$-form of the map $\pr_1$ from above. Consequently $\sh^2(\pi_1\big(\rR_{L/K}(X)/K\big)$ equals $\pi_1(X/L)$ and we have the following structure theorem for the fundamental group of a Weil restriction of scalars.
\begin{thm}  \label{thm:pi1weiltext}
The fundamental group $\pi_1(\rR_{L/K} X/K)$ is isomorphic to the non-abelian induction $\Ind_{\Gal_L}^{\Gal_K} \pi_1(X/L)$.
\end{thm}


\end{document}